\documentclass[11pt,amsfonts, epsfig]{amsart}
\usepackage{amsmath, amscd, amssymb}
\usepackage{graphicx, psfrag} 
\usepackage{graphpap, color}
\usepackage{mathrsfs}
\usepackage{pstricks}
\usepackage{color}
\usepackage{cancel}
\usepackage[mathscr]{eucal}
\usepackage{pstricks}
\usepackage{cancel}
\usepackage{verbatim}
\usepackage{latexsym}
\usepackage[all]{xy}

\usepackage{lipsum}
\makeatletter
\g@addto@macro{\endabstract}{\@setabstract}
\newcommand{\authorfootnotes}{\renewcommand\thefootnote{\@fnsymbol\c@footnote}}%
\makeatother

\usepackage{upgreek}

 \usepackage{wasysym}
\usepackage{ esint }

\numberwithin{equation}{section}

\def\sO{{\mathscr O}}
\def\sC{{\mathscr C}}

\def\sN{{\mathscr N}}
\def\sL{{\mathscr L}}

\def\sO{\mathscr{O}}

\newcommand{\PP}{\mathbb{P}}

\newcommand{\cal}{\mathcal}

\def\cL{{\cal L}}

\def\v1{{\vec{1}}}

\def\begeq{\begin{equation}}
\def\endeq{\end{equation}}
\def\and{\quad{\rm and}\quad}

\def\sub{\subset}

\def\and{\quad\text{and}\quad}

\newtheorem{prop}{Proposition}[section]

\newtheorem{theo}[prop]{Theorem}

\newtheorem{rema}[prop]{Remark}
\newtheorem{exam}[prop]{Example}

\newtheorem{conj}[prop]{Conjecture}

\newtheorem{defi-prop}[prop]{Definition-Proposition}
\newtheorem{defi-theo}[prop]{Definition-Theorem}

\def\sS{\mathscr S}

\def\sO{{\mathscr O}}

\def\sD{{\mathscr D}}

\def\beq{\begin{equation}}
\def\eeq{\end{equation}}

\def\Pf{{\PP^4}}

\def\bee{\begin{equation}}
\def\eeq{\end{equation}}

\def\sC{{\mathscr C}}

\def\bd{{\mathbf d}}

\def\bmu{{\boldsymbol \mu}}

\begin{document}

\title[MSP]{A Survey on Mixed Spin P-Fields}{
\author[Huai-Liang Chang]{Huai-Liang Chang$^1$}
\address{Department of Mathematics, Hong Kong University of Science and Technology, Hong Kong} \email{mahlchang@ust.hk}
\thanks{${}^1$Partially supported by  Hong Kong GRF grant 600711 and 6301515}

\author[Jun Li]{Jun Li$^2$}
\address{Shanghai Center for Mathematical Sciences, Fudan University, China; \hfil\newline 
\indent Department of Mathematics, Stanford University,
USA} \email{jli@math.stanford.edu}
\thanks{${}^2$Partially supported by  NSF grant 
NSF-1104553 and DMS-1159156}

\author[Wei-Ping Li]{Wei-Ping Li$^3$}
\address{Department of Mathematics, Hong Kong University of Science and Technology, Hong Kong} \email{mawpli@ust.hk}
\thanks{${}^3$Partially supported by by Hong Kong GRF grant 602512 and 6301515}

\author[Chiu-Chu Melissa Liu]{Chiu-Chu Melissa Liu$^4$}
\address{Mathematics Department, Columbia University}
 \email{ccliu@math.columbia.edu}
\thanks{${}^4$Partially supported by  NSF grant DMS-1206667 and DMS-1159416
}

\maketitle

\section{Gromov-Witten invariants of quintics}

For the Fermat quintic polynomial $W_5=x_1^5+\ldots+x_5^5$,  the counting of genus $g$ curves of   degree $d$ on the  quintic Calabi-Yau three-fold $Q=\{W_5=0\}\subset \mathbb P^4$ is a challenging problem in enumerative geometry. Since the seminal paper of Candelas, dela Ossa, Green and Parkes \cite{Can},
a modified version has been intensively studied in string theory  as well as algebraic geometry via the stable maps of Kontsevich and virtual cycles theory developed by Li-Tian \cite{LT} and Behrend-Fantachi \cite{BF}.

For $d, g\in \mathbb Z$, the moduli space of stable maps from genus $g$ nodal curves to $Q$ of degree $d$ is 
$${
\overline{M}_{g}(X,d)=\{[f:C\to X]\mid C\ \text{nodal}, \,g(C)=g,\, f_*[C]=d,\, \text{Aut}  (f)<\infty\}.}
$$
The Gromov-Witten invariants are defined as
$$N_{g,d}\colon =\int_{[\overline{M}_g(X, d)]^{vir}}1\in\mathbb Q.$$
One of the main unsolved problems in Gromov-Witten theory is to determine 
$$F_g(q)\colon=\sum_{d}N_{g,d}q^d.$$

From the Super-String Theory side, in 1991 Candelas et.al. found a closed formula for genus zero $F_0(q)$ using $T$-duality and mirror symmetry (\cite{Can}). In 1993, Bershadsky, Cecotti, Ooguri and Vafa developed the Kodaira-Spencer theory and determined the genus one $F_1(q)$ (\cite{BCOV}). For higher genus, in 2009 Huang, Klemm and Quackenbush determined $F_g(q)$ for $g$ up to $51$ (\cite{HKQ}).

From the mathematical side, Kontsevich derived a torus localization to calculate the genus zero GW-invariants $N_{0,d}$. Givental \cite{Gi}, Lian-Liu-Yau \cite{LLY} determined the genus zero case $F_0(q)$. Later on, more people  worked on this topic. The genus one case $F_1(q)$ was solved in 2000's. The second named author and Zinger in \cite{LZ} obtained a formula $N_{1,d}^{red}=N_{1,d}-\displaystyle{\frac{1}{12}}N_{0,d}$ where $N_{1,d}^{red}$ is certain reduced GW-invariants. Using this formula and $\mathbb C^*$-localization, Zinger in \cite{Zi} succeeded in determining $F_1(q)$. Gathmann \cite {Gath} provided an algorithm for $N_{1,d}$ using the relative GW-invariant formula. For higher genus case, Maulik and Pandharipande {\black found an algorithm \cite[Section 3.2]{MP}} using the { algebraic version of degeneration formula
{\black \cite{Lideg}} (see analogue formula \cite{LR})} and used it for some theoretical applications.
 {Despite these progress, a lot of questions on higher genus GW invariants of quintic Calabi-Yau threefolds remains open.}

{ It remains a central problem in Gromov-Witten theory  {to develop new techniques} to
calculate all genus GW-invariants of quintic Calabi-Yau threefolds.}


\section{ Witten's vision and FJRW invariants}

\subsection{Witten's vision}

The same quintic polynomial $W_5=x_1^5+\ldots+x_5^5$ can also give a map $\mathbb C^5\to \mathbb C$. The corresponding physical theory is the Landau-Ginzburg theory. In \cite{GLSM}, Witten studied phase transitions involving GW theory on the quintic $Q$ and the LG-model for $W_5$. Mathematically, the set-up is as follows. Let $\mathbb C^*$ act on 
$$\mathbb C^6=\mathbb C^5\times \mathbb C=\{(x_1, \ldots, x_5, p)\}
$$ with weights $(1, \ldots, 1, -5)$. Then the map $p\cdot (x_1^5+\ldots+x_5^5)\colon \mathbb C^6\to\mathbb C$ is $\mathbb C^*$-equivariant. The quotient $[\mathbb C^6/\mathbb C^*]$ has two GIT quotients:
$$
\big((\mathbb C^5-\{\vec 0\})\times \mathbb C\big)/\mathbb C^*=K_{\mathbb P^4},
$$
and
$$
\big(\mathbb C^5\times (\mathbb C-0)\big)/\mathbb C^*=[\mathbb C^5/\mathbb Z_5].
$$ 
Here $[\mathbb C^5/\mathbb Z_5]$ represents the quotient stack. 
The field theory valued in $K_{\mathbb P^4}$ is the GW theory of the quintic $Q$ and the field theory valued in $[\mathbb C^5/\mathbb Z_5]$ leads to the Witten's spin class. The latter was generalized to quasi-homogeneous polynomials by Fan, Jarvis and Ruan \cite{FJR1, FJR2}, which is called FJRW theory. Witten's vision is that these two theories are related via 
a phase transition. 

\subsection{P-fields treatment of GW and FJRW}

The notion of P-fields was introduced by Guffin and Sharpe in \cite{Sharpe} for genus zero LG-theory of $(K_{\mathbb P^4}, W_5)$. Mathematically, The first and second named authors developed the theory of P-fields for all genus GW invariants.

We start with LG-theory for $K_{\mathbb P^4}$.
A field taking values in $\big((\mathbb C^5-\{\vec 0\})\times \mathbb C\big)/\mathbb C^*$ is 
$$\xi=(\sC, \sL, \varphi_1, \ldots, \varphi_5, \rho)$$
where $\sC$ is a complete nodal curve, $\sL$ is an invertible sheaf on $\sC$,  $\varphi_i\in H^0(\sC, \sL)$, and $\rho\in H^0(\sL^{\vee5}\otimes \omega_\sC)$. Since the weights of the action $\mathbb C^*$ on $\mathbb C^5$ and $\mathbb C$ are $(1, \ldots, 1)$ and $-5$ respectively, while $\varphi_i$ is a section of $\sL$,  a priori $\rho$ has to be a section of $\sL^{\vee 5}$ , but we choose $\rho$ to be a section  of $\sL^{\vee 5}\otimes \omega_{\sC}$  due to a technical reason. Since we deleted the origin $\vec 0$  from $\mathbb C^5$, $(\varphi_1, \ldots, \varphi_5)$ must be nowhere zero. Since we consider the quotient space via the  $\mathbb C^*$-action, we need to introduce the natural equivalence $(\sC, \sL)\to (\sC, \sL)$ where { $\mathbb C^*$} acts on $\sL$ via scalar multiplication. Finally we say $\xi$ is stable if Aut($\xi$) is finite.

In fact, what we have gotten so far  is a stable map to $\mathbb P^4$ with a P-field. The moduli space of such objects is
\begin{eqnarray*}
\overline{M}_g(\mathbf P^4,d)^p=\{[f,\sC,\rho]\mid [f,\sC]\in \overline{M}_g(\mathbf P^4,d), \,
\rho \in H^0(\sC,f^*\mathcal O(5)\otimes\omega_\sC)\}.
\end{eqnarray*}
Note that the data $([f,\sC],\rho)$ is equivalent to the data $ (\sC,\cL,\varphi_1,\cdots,\varphi_5, \rho)$ since  the map 
$f $ is equivalent to the line bundle $\sL=f^*\sO_{\mathbb P^4}(1)$ with five sections $(\varphi_1,\cdots,\varphi_5)$ of $\sL$.

The first and second named authors constructed the GW invariants of stable maps with P-fields as follows. The moduli stack $\overline{M}_g(\mathbf P^4,d)^p$, relative to the stack $\sD=\{(\sC, \sL)\}$,  has a perfect obstruction theory. At $\xi=(\sC, \sL, \varphi_i,\rho)$, the obstruction sheaf restricted to $\xi$ is 
$$\sO b|_{\xi}=H^1(\sL)^{\oplus 5}\oplus H^1(\sL^{\vee 5}\otimes \omega_\sC).$$
There exists a cosection 
$$\sigma\colon \sO b\to \sO_{\overline{M}_g(\mathbf P^4,d)^p}$$
constructed as follows. Let 
$$(\dot\varphi_1,\ldots, \dot\varphi_5, \dot\rho)\in H^1(\sL)^{\oplus 5}\oplus H^1(\sL^{\vee 5}\otimes \omega_\sC)=\sO b|_{\xi}.
$$

Define  $\sigma|_{\xi}(\dot\varphi_1,\ldots, \dot\varphi_5,\dot\rho)\colon =\dot\rho\sum_{i=1}^5\varphi_i^5+\rho\sum_{i=1}^55\varphi_i^4\dot\varphi_i$.

The degeneracy locus $D(\sigma)$ of the cosection consists of $\xi$ such that $\sigma|_{\xi}$ is zero, i.e., $\sigma|_{\xi}(\dot\varphi_1,\ldots, \dot\varphi_5,\dot\rho)=0$ for all $\dot\varphi_i$ and $\dot\rho$. Thus 
$$D(\sigma)=\{\xi \in \overline{M}_g(\mathbf P^4,d)^p\, |\, \rho=0\hbox{ and } \sum_{i=1}^5\varphi_i^5=0\}=\overline{M}_g(Q,d)\subset \overline {M}_g(\mathbf P^4,d).$$
The expression of the cosection $\sigma$ comes from taking  the derivative of $p\cdot W_5=p(x_1^5+\ldots+x_5^5)$ with respect to  the time variable $t$ where $p$ and $x_i$ are regarded as functions of $t$ following physical notations. 
 
Since  $\rho$ is a section, the moduli space $\overline {M}_g(\mathbf P^4,d)^p$ is not proper (when $g\ge 1$) and hence cannot be used to define invariants. However, {Kiem and the second named author
\cite{KL} developed a theory of cosection localization virtual cycles which, applied to this case, resolves the non-proper
issue. More precisely, one checks that the degeneracy locus $D(\sigma)$ is the moduli space of stable maps to the quintic 
$Q$ and thus proper.} 
 \begin{theo}[H.L. Chang - J. Li \cite{CL}]   Using the cosection localized virtual cycle, one obtains the cycle
 $$[ \overline{M}_g(\mathbf P^4,d)^p]^{vir}_{loc}\in A_*D(\sigma)=A_*\overline{M}_g(Q,d).
 $$
Furthermore, let 
P-fields GW invariants $N^p_{g, d}=\int_{[\overline{M}_g(\mathbf P^4,d)^p]^{vir}_{loc}}1\in \mathbb Q$, then
 $$N_{g,d}=(-1)^{d+g+1}N^p_{g, d}.$$
 \end{theo}
 
 The advantage of this result is that $F_g(q)=\sum_d N_{g, d}q^d$ now becomes a topological string  amplitude of a field theory valued in $K_{\mathbb P^4}=\big((\mathbb C^5-\vec 0)\times \mathbb C\big)/\mathbb C^*$.

Now let's consider the field theory valued in $[\mathbb C^5/\mathbb Z_5]$. This theory originated from Witten's spin class \cite{GLSM}. Its algebraic constructions (in narrow case) were given by Polishchuk-Vaintrob \cite{PV} and Chiodo  \cite{Chi}. The full theory was developed by Fan, Jarvis and Ruan \cite{FJR1, FJR2}, known as FJRW theory. We will touch (narrow) FJRW invariants following the construction by the first, second and third named authors \cite{CLL}.

As the case in the P-fields treatment of GW theory, a field in $\big(\mathbb C^5\times (\mathbb C-0)\big)/\mathbb C^*=[\mathbb C^5/\mathbb Z_5]$ consists of 
$$\xi=(\Sigma^{\sC}, \sC, \sL, \varphi_1,\ldots,\varphi_5, \rho)$$
where $(\Sigma^{\sC},\sC)$ is a pointed twisted curve with markings $\Sigma^{\sC}$ possibly stacky, $\sL$ is an invertible sheaf on $\sC$, $\varphi_i\in H^0(\sL)$, and $\rho\in H^0(\sL^{\vee 5}\otimes \omega^{\log}_{\sC})$ with $\omega^{\log}_{\sC}=\omega_{\sC}(\Sigma^{\sC})$. Since we deleted the origin in $\mathbb C$, the section $\rho$ must be nowhere vanishing and hence $\sL^{\vee 5}\otimes \omega_\sC\cong \sO_\sC$, or equivalently $\sL^{\otimes 5}\cong \omega^{\log}_{\sC}$. Therefore $(\Sigma^{\sC}, \sC, \sL)$ is a $5$-spin curve. $(\varphi_1, \ldots, \varphi_5)$ give five fields. Thus we get a moduli space of $5$-spin curves with five fields:
\begin{eqnarray*}
\overline{M}_{g,\gamma}^{1/5,5p}=\{(\sC, \Sigma^{\sC}, \cL, \varphi_1,\cdots,\varphi_5, \rho) \mid \text{$\rho$ is nowhere zero}\}.
\end{eqnarray*}
Here $\gamma$ is the monodromy data: if $\Sigma_1$ is a stacky marking on $\sC$, then $\bmu_5$ acts on $\sL|_{\Sigma_1}$ with weight 
$\gamma_1={ \exp}(2\pi i r/5)$ where $0\le r\le 4$. Narrow means $0<r\le 4$. If $\Sigma_1$ is an ordinary marking, $\gamma_1$ is taken to be $1$. 

Similar to GW case, the moduli stack $\overline{M}_{g,\gamma}^{1/5,5p}$, relative to the stack $\sD=\{(\Sigma^{\sC}, \sC, \sL)\}$, has a perfect obstruction theory. There exists a cosection
$\sigma\colon \sO b\to \sO_{\overline{M}_{g,\gamma}^{1/5,5p}}$ whose degeneracy locus  is 
$$D(\sigma)=\{\xi\in \overline{M}_{g,\gamma}^{1/5,5p} \,|\, \varphi_i=0 \hbox{ for all $i$}\}=\overline{M}_{g,\gamma}^{1/5}=\{(\Sigma^{\sC}, \sC, \sL)\, |\, \sL^{\otimes 5}\cong \omega^{\log}_\sC\},$$
which is the moduli space of $5$-spin curves. 
\begin{theo}[H.L. Chang - J. Li - W.P. Li \cite{CLL}] The (narrow) FJRW invariants can be constructed using cosection localized virtual cycles of the moduli space of spin curves with five fields:
$$[\overline{M}_{g,\gamma}^{1/5,5p}]^{vir}_{loc}\in A_*\overline{M}_{g,\gamma}^{1/5}.$$
\end{theo}

This construction is an algebraic geometric version of Witten's original construction. Witten considered the moduli space of $5$-spin curves $(\Sigma^{\sC}, \sC, \sL)$ with smooth sections. For our set-up, the corresponding 
  Witten equations are given by 
\begin{eqnarray}\label{Witten's-eq}
\bar\partial s_i+\overline{\partial_{x_i}W_5(s_1, \ldots, s_5)}=0, \quad i.e., \quad \bar\partial s_i+5\overline {s_i^4}=0.
\end{eqnarray}
This is used to construct Witten's top Chern class to define invariants on the moduli space of $5$-spin curves. From Witten's equation, the term $\bar \partial s_i$ gives the obstruction class to extend a holomorphic section. Thus the left hand side of (\ref{Witten's-eq}) gives a (differential) section of the obstruction sheaf of the moduli of spin curves with fields. Now substitute the complex conjugate in the Witten's equation by the Serre duality, the LHS of (\ref{Witten's-eq})  becomes the cosection. 

There is an important subclass of FJRW invariants: those with the  insertion $-\displaystyle\frac{2}{5}$. Let $\sC$ have $k$ markings with all $\gamma_j=\zeta^2$ for $1\le j\le k$ where 
$\zeta={ \exp}(2\pi i /5)$. Define
$$
\Theta_{g, k}\colon =\int_{[\overline{M}_{g,(\gamma_j)}^{1/5,5p}]^{vir}_{loc}}1\in\mathbb Q, \quad \hbox{for $k+2-2g=0$ {mod} $5$}.
$$

{\black  It is shown \cite{CLLL2} that $\{\Theta_{g,k}\}_{g,k}$ determine all FJRW invariants with descendents (for the quintic singularity), where an explicit formula will be given in \cite{twFJRW}. For this reason we call $\{\Theta_{g,k}\}_{g,k}$ the primary FJRW invariants. }

\section{Master space technique and mixed spin fields}

In the previous section, we discussed the LG-field theoretic description of GW theory of the quintic and FJRW theory of $(\mathbb C^5, W_5)$. Witten's vision is to link these two theories via { a phase transition} with respect to some complexified parameter. {The approach by the authors} is to develop a field theory valued in the master space to geometrically realize the ``wall-crossings" of these two field theories. 

\subsection{Master space technique} Now we explain the master space technique to understand the wall-crossings between $K_{\mathbb P^4}$ and $[\mathbb C^5/\mathbb Z_5]$. 

Consider a $\mathbb C^*$-action on $\mathbb C^5\times\mathbb C\times \mathbb P^1$, for $t\in \mathbb C^*$, 
$$(x_1, \ldots, x_5, p,[u_1, u_2])^t\colon =(tx_1,\ldots, tx_5, t^{-5}p, [tu_1, u_2]).$$
It has a GIT quotient 
$$
W\colon =(\mathbb C^5\times \mathbb  C\times \mathbb P^1-\sS)/\mathbb C^*\hbox{ where }  \sS\colon=\{(x_i=0=u_1)\cup (\rho=0=u_2)\}.
$$

Consider a $\mathbb C^*$-action on $W$ and, to avoid confusions, we call this action $T$-action. For $t\in T=\mathbb C^*$,
$$
(x_1,\ldots, x_5, p, [u_1, u_2])^t=(x_1,\ldots, x_5, p, [tu_1, u_2]).
$$
The $T$-fixed locus is
$$
W^T=K_{\mathbb P^4}{ \times} \{0\}\coprod \vec 0\times 
{ \left( (\mathbb P^1-\{0,\infty\})/{\mathbb C}^*\right)  } \coprod 
[\mathbb C^5/\mathbb Z_5]\times \{\infty\}
$$
where $0=[1, 0]$ and $\infty=[1, 0]$ in $\mathbb P^1$.

Take a $T$-equivariant form $\mu$ on $W$, then we have
$$
0=\left[\int_W\mu\cap c_1({\bf 1}_{wt=1})\right]_0=-\int_{K_{\mathbb P^4}}\mu+\int_{[\mathbb C^5/\mathbb Z_5]}\mu
+\left[\int_{ \{ \text{point} \} }\frac{\mu}{*}\right]_0
$$
where ${\bf 1}_{wt=1}$ is the $T$-linearized trivial line bundle with weight $1$ and $[\ldots]_0$ means taking degree zero part in the equivariant parameter. Thus the wall-crossing can be expressed as
$$
\int_{[\mathbb C^5/\mathbb Z_5]}\mu -\int_{K_{\mathbb P^4}}\mu=\text{error} 
=\left[\int_{ \{ \text{point} \} }\frac{\mu}{*}\right]_0$$

\subsection{Mixed spin P-fields}
Now we consider a field theory valued in $W$. Similar to the case of the field theory of GW valued in $K_{\mathbb P^4}$, 
{the authors introduced the notion of mixed spin $P$-fields (MSP for short) (\cite{CLLL}). An MSP field is}
$$\xi=(\Sigma^\sC, \sC, \sL , \sN, \varphi_1,\ldots, \varphi_5, \rho, \nu=[\nu_1,\nu_2]).$$
$(\Sigma^{\sC}, \sC)$ is a pointed twisted curve. $\sL$ and $\sN$ are invertible sheaves on $\sC$. $\sL$ is as before but $\sN$ is new due to the extra factor $\mathbb P^1$ in the master space technique. $\varphi_i\in H^0(\sL)$ and $\rho\in H^0(\sL^{\vee 5}\otimes \omega_{\sC}^{\log})$ as before. $\nu_1\in H^0(\sL\otimes \sN)$ and $\nu_2\in H^0(\sN)$. $\nu=[\nu_1, \nu_2]$ is a new field. We also have a narrow condition: $\varphi_i|_{\Sigma^{\sC}}=0$. There are combined GIT-like stability requirements: $(\varphi_1, \ldots, \varphi_5, \nu_1)$ is nowhere vanishing coming from excluding 
$\{(x_i=0=u_1)\}$ in $W$; $(\rho, \nu_2)$ is nowhere vanishing coming from excluding $\{(\rho=0=u_2)\}$ in $W$; and $(\nu_1, \nu_2)$ is nowhere vanishing coming for $[u_1, u_2]\in \mathbb P^1$. We say $\xi$ is stable if Aut($\xi$)  is finite. For simplicity, we use $\varphi$ to represent $(\varphi_1, \ldots, \varphi_5)$.

In order to understand why the moduli space of MSP fields geometrically contains the moduli space of stable maps with P-fields and the moduli space of spin curves with five P-fields, we examine the moduli space of MSP fields in details. 

Let $\xi$ be a MSP field. When $\nu_1=0$, since $(\varphi_1, \ldots, \varphi_5,\nu_1)$ is nowhere zero, we must have that $(\varphi_1, \ldots, \varphi_5)$ is nowhere zero. Since $(\nu_1, \nu_2)$ is nowhere zero, $\nu_2$ must be nowhere zero. Since $\nu_2$ is a section of $\sN$, $\sN\cong \sO_\sC$. There is no restriction on $\rho$. Thus $\xi\in \overline{M}_g(\mathbf P^4,d)^p$ and we get GW theory of the quintic $Q$.

When $\nu_2=0$, since $(\rho, \nu_2)$ is nowhere zero, $\rho$ must be nowhere vanishing. Since $\rho$ is a section of $\sL^{\vee 5}\otimes \omega_\sC^{\log}$, we must have $\sL^5\cong \omega_\sC^{\log}$. Also $\nu_1$ must be nowhere zero. Thus $\sL\otimes \sN\cong\sO_\sC$, i.e., $\sN\cong \sL^\vee$.  $\varphi_1,\ldots, \varphi_5$ can be arbitrary. Thus $\xi\in \overline{M}_{g,(\gamma_j)}^{1/5,5p}$ and we get FJRW theory. 

When $\rho=0$ and $\varphi_i=0$ for $1\le i\le 5$, $\nu_1, \nu_2$ must be nowhere zero. Thus $\sN\cong \sO_\sC$ and $\sL\cong\sO_\sC$. Hence we get stable curves. 

\begin{theo}[
\cite{CLLL}] The moduli stack $W_{g, \gamma, {\bf d}}$ of stable MSP fields of genus $g$, monodromy $\gamma=(\gamma_1, \ldots,\gamma_\ell)$ of $\sL$ along $\Sigma^{\sC}$, and degree ${\bf d}=(d_0, d_\infty)$ of $\sL\otimes \sN$ and $\sN$ respectively,  is a separated DM stack of locally finite type.
\end{theo}

The moduli stack $W_{g, \gamma, {\bf d}}$ admits a natural $\mathbb C^*$-action: for $t\in \mathbb C^*$, 
\begin{eqnarray*}
(\Sigma^{\sC}, \sC, \sL, \sN, \varphi, \rho, \nu_1, \nu_2)^t\colon =(\Sigma^{\sC}, \sC, \sL, \sN, \varphi, \rho, t\nu_1, \nu_2).
\end{eqnarray*}
$W_{g, \gamma, {\bf d}}$ is not  proper since $\varphi$ and $\rho$ are sections of invertible sheaves. Thus we cannot do integrations on this stack. However, there exists a cosection of its obstruction sheaf. Using the arguments similar to  GW case and LG case, we have the following theorem.

\begin{theo}[\cite{CLLL}] The moduli stack $W_{g, \gamma,  {\bf d}}$ has a $\mathbb C^*$-equivariant perfect obstruction theory, an equivariant cosection $\sigma$ of its obstruction sheaf, and thus carries an equivariant cosection localized virtual cycle
$$[W_{g, \gamma, {\bf d}}]^{vir}_{loc}\in A_*^{\mathbb C^*} W_{g,\gamma,{\bf d}}^-$$
where $W_{g,\gamma,{\bf d}}^-$ is the degeneracy locus of $\sigma$, i.e.,
$$W_{g,\gamma,{\bf d}}^-\colon=(\sigma=0)=\{\xi\in W_{g,\gamma,{\bf d}}\, |\, \sC=(\varphi=0)\cup (\varphi_1^5+\ldots+\varphi_5^5=0=\rho)\}.$$
\end{theo}

In order to do integration on $W_{g,\gamma,{\bf d}}^-$, one needs it to be proper. In fact, we have
\begin{theo}[\cite{CLLL}] The degeneracy locus $W_{g,\gamma,{\bf d}}^-$ is a proper $\mathbb C^*$-DM stack of finite type.
\end{theo}

{\black From the proof of the  properness,  we see  a phenomenon which creates line bundles' spin structures in the LG-phase via a limit of a family of P-fields in CY-phase. We give this phenomenon the name ``Landau-Ginzburg transition"(or CY-to-LG transition). It is under this phenomenon that  FJRW theory captures the {ghosts' contributions in  GW theory (\footnote{ A map $f$ from a curve $C$ to $\Pf$ or $Q$ is called a ghost if there are positive-genus components of $C$ that are contracted to points by $f$. Over ghosts,   P-fields can be nonvanishing,  and such P-fields contribute to GW invariants of the quintic  as ``counting ghosts". For example,  in Li-Zinger formula
$N_{1,d}=N_{1,d}^{red}+\displaystyle{\frac{1}{12}}N_{0,d}$, 
the number $1/12$ comes from the contribution by P-fields. When genus increases,  such contribution is difficult to locate. 
MSP program provides a platform that ghosts' contributions can be captured in another phase (LG-phase)  instead.\black
}) }in the realm of MSP moduli.}

 \begin{exam}\label{g1d1} The graph of fixed  points of $W_{1, \emptyset, (1, 0)}^-$. 
\end{exam}
 
  So the curves are  elliptic curves without markings, $\deg\sL=1$, and $\deg \sN=0$. The graph type of fixed points which have contributions to the computations are of the following  four types { $\Gamma_1$, 
$\Gamma_2$, $\Gamma_3$, $\Gamma_4$}:

\begin{figure}[h] 
\begin{center}
\psfrag{G1}{$\Gamma_1$}
\psfrag{G2}{$\Gamma_2$}
\psfrag{G3}{$\Gamma_3$}
\psfrag{G4}{$\Gamma_4$}
\psfrag{0}{\small $0$}
\psfrag{1}{\small $1$}
\psfrag{infty}{\small $\infty$}
\includegraphics[scale=0.6]{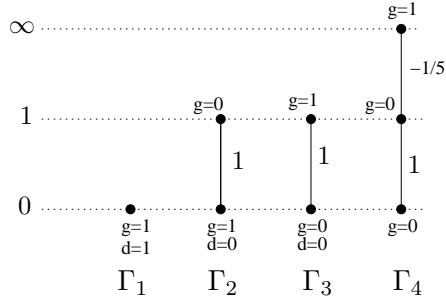}
\end{center}
\caption{graphs for $g=1$, $\gamma=\emptyset$, $\bd=(1,0)$.}
\end{figure}
  
{Some explanations of the figure are in order.} 
The bottom horizontal line corresponds to $\nu_1=0$. The middle horizontal line corresponds to $\varphi=0$ and $\rho=0$. The top horizontal line corresponds to $\nu_2=0$. A vertex represents a connected curve if it is stable (called a stable vertex), or a node if it has two edges, or nothing if it has only one edge attached to it. For each vertex, $g$ represents the genus of the curve. If it is not a stable vertex, nor a node, we use $g=0$ here even though it doesn't represent a rational curve.  An edge is a rational curve. The number near an edge is the degree of $\sL$ on the curve. A stable vertex on the horizontal lines $0, 1, \text{infty}$ means $\nu_1=0$, $\varphi=0=\rho$, or $\nu_2=0$ on the curve respectively. 

To be more precise, the graph { $\Gamma_1$} represents an elliptic curve with degree 1 line bundle $\sL$ on the curve and $\nu_1=0$ on the whole curve. So it corresponds to stable maps from elliptic curves to the quintic $Q$ with degree $1$. { The graph  $\Gamma_2$} represents a union of an elliptic curve $E$ with a rational curve $E_0$ intersecting at one node. $E$ is a stable vertex on the bottom horizontal line. $E_0$ is the edge. $\deg\sL$ is $0$ on the elliptic curve and $1$ on the rational curve. 
{ The graph $\Gamma_3$}  is similar to $\Gamma_2$, a union of an elliptic curve $E$ and a rational curve $E_0$ intersecting at one node. $E$ is a stable vertex on the middle horizontal line and $E_0$ is the edge. On $E$, $\varphi=0=\rho$ and both $\sL$ and $\sN$ are trivial. Thus it represents the moduli space of elliptic curves with one marking coming from the node. { The graph $\Gamma_4$} represents a union of two rational curves $E_0$ and $E_\infty$ and an elliptic curve $E$. On the lower edge $E_0$, $\deg(\sL|_{E_0})=1$ and $\rho|_{E_0}=0$. Note that on each irreducible component, either $\rho=0$ or $\varphi=0$. On $E_\infty$, $\deg(\sL|_{E_\infty})=-1/5$ and $\varphi|_{E_\infty}=0$. $E_\infty$ and $E_0$ intersect at one node. $E_\infty$ is a twisted curve intersecting the elliptic curve $E$ at a stacky point. Thus $\sL|_{E_\infty}$ is an invertible sheaf on the twisted curve $E_\infty$. $E$ is a spin elliptic curve with one marking from the node. Thus $\deg(\sL|_E)=1/5$.

\begin{figure}[h] 
\begin{center}
\psfrag{G1}{\footnotesize $\Gamma_1$}
\psfrag{G2}{\footnotesize $\Gamma_2$}
\psfrag{G3}{\footnotesize $\Gamma_3$}
\psfrag{G4}{\footnotesize $\Gamma_4$}
\psfrag{G5}{\footnotesize $\Gamma_5$}
\psfrag{G6}{\footnotesize $\Gamma_6$}
\psfrag{G7}{\footnotesize $\Gamma_7$}
\psfrag{G8}{\footnotesize $\Gamma_8$}
\psfrag{G9}{\footnotesize $\Gamma_9$}
\psfrag{G10}{\footnotesize $\Gamma_{10}$}
\psfrag{G11}{\footnotesize $\Gamma_{11}$}
\psfrag{G12}{\footnotesize $\Gamma_{12}$}
\psfrag{G13}{\footnotesize $\Gamma_{13}$}
\psfrag{G14}{\footnotesize $\Gamma_{14}$}
\psfrag{G15}{\footnotesize $\Gamma_{15}$}
\psfrag{0}{\footnotesize $0$}
\psfrag{01}{\footnotesize $1$}
\psfrag{infty}{\footnotesize $\infty$}
\includegraphics[scale=0.54]{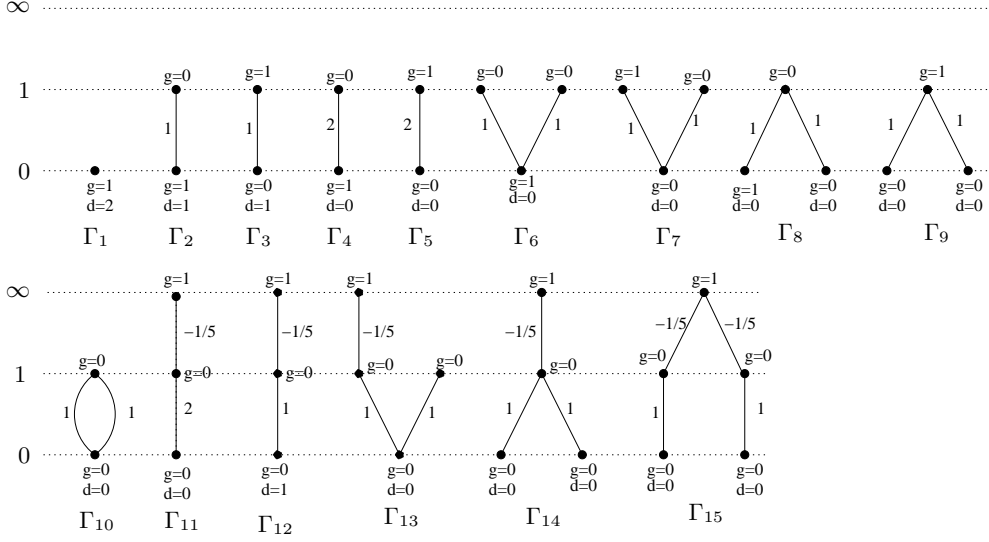}
\end{center}
\caption{graphs for $g=1$, $\gamma=\emptyset$, $\bd=(2,0)$.}
\end{figure}

\begin{exam}\label{g1d2} The graph of fixed  points of $W_{1, \emptyset, (2, 0)}^-$. 
\end{exam}

{ In this case} the curve is an elliptic curve without markings, $\deg\sL=2$ and $\deg \sN=0$. The graph type of fixed points which have contributions to the computations have15 types, listed from { $\Gamma_1$} to 
{ $\Gamma_{15}$}.

\section{Vanishing and polynomial relations}

How to extract information of GW  and/or FJRW invariants from the cycle $[W_{g, \gamma, \bf d}]^{vir}_{loc}$? Let's consider a less general case to illustrate  key ideas. Take $\gamma=\emptyset$, i.e., no markings. Then by virtual dimension counting, we have
$$[W_{g, \bf d}]^{vir}_{loc}\in H^{\mathbb C^*}_{2(d_0+d_\infty+1-g)}(W_{g, \bf d}^-, \mathbb Q).$$
When $d_0+d_{\infty}+1-g>0$, letting $u=c_1({\bf 1}|_{wt=1})$, i.e. $u$ is the parameter for $H_{\mathbb C^*}^*(pt)$, we have
\begin{eqnarray*}
[u^{d_0+d_\infty+1-g}\cdot [W_{g, \bf d}]^{vir}_{loc}]_0=0.
\end{eqnarray*}
Here $[\cdot]_0$ is the degree zero term in the variable $u$. 

Let $\Gamma$ be a graph associated to fixed points of the $\mathbb C^*$-action of $W_{g, \bf d}$ and $F_{\Gamma}$ be a connected component of $W_{g, \bf d}^{\mathbb C^*}$ of the graph type $\Gamma$. Apply the cosection localized version proved in \cite{CKL} of the virtual localization formula 
in \cite{GP}, we have
\begin{eqnarray}\label{virtual-local}
\sum_{\Gamma}\left[u^{d_0+d_\infty+1-g}\frac{[F_\Gamma]^{vir}_{loc}}{e(N_{F_\Gamma})}\right]_0=0.
\end{eqnarray}
 To deal with $[F_\Gamma]^{vir}_{loc}$, we need a decomposition result to be explained below.

Let, again, $\xi=(\sC,  \sL, \sN, \varphi, \rho, \nu_1,\nu_2) \in (W_{g, d})^{\mathbb C^*}$ be a MSP field fixed by the $\mathbb C^*$-action. We set
\begin{enumerate}

\item $\sC_0$ to be the part of $\sC$ where $\nu_1=0$;

\item $\sC_1$ to be the part of $\sC$ where $\varphi=0=\rho$ and hence $\nu_1=1 =\nu_2$, i.e., $\nu_1$ and $\nu_2$ are nowhere zero;

\item $\sC_\infty$ to be the part of $\sC$ where $\nu_2=0$.

\end{enumerate}

Thus $\xi|_{\hbox{(connected component of $\sC_0$})}$ is in $\overline{M}_{g^\prime, n^\prime}(\mathbb P^4, d^\prime)^p$ which gives Gromov-Witten  invariants.  Here marked points appear coming from some nodes on $\sC_0$. $\xi|_{\hbox{(connected component of $\sC_1$})}$ is in $\overline{M}_{g^\prime, n^\prime}$ which gives Hodge integrals. $\xi|_{\hbox{(connected component of $\sC_\infty$})}$ is in $\overline{M}^{\frac{1}{5}, 5p}_{g^\prime, \gamma^\prime}$ which gives FJRW invariants where $\gamma^\prime$ appears because of some stacky nodes on $\sC_\infty$. 

We have the following decomposition result:
\begin{eqnarray*}
[F_\Gamma]^{vir}_{loc}=c\prod [\hbox{moduli of $\xi|_{\sC_0}$]}^{vir}_{loc}\cdot[\hbox{moduli of $\xi|_{\sC_1}$}]^{vir}_{loc}\cdot [ \hbox{moduli of $\xi|_{\sC_\infty}$}]^{vir}_{loc}
\end{eqnarray*}
where $c$ is a constant.  The first factor gives GW invariants of stable maps to $\PP^4$ with $P$-fields, i.e. $N_{g',d'}$. The second  factor gives Hodge integrals on $\overline M_{g',n'}$. The third factor gives FJRW invariants of insertions $-\frac{2}{5}$ (after using a vanishing).  After $e(N_{F_\Gamma})$'s are calculated, using the polynomial relations (\ref{virtual-local}), we obtain the following results about GW invariants of the quintic. 
\begin{theo}[\cite{CLLL2}]\label{thm-induction}
Letting $d_\infty=0$, the relations (\ref{virtual-local}) provide an effective algorithm to evaluate GW invariants $N_{g,d}$
provided the following are known
\begin{enumerate}
\item $N_{g',d'}$ for $(g',d')$ such that $g'< g$, and $d'\le d$;
\item $N_{g,d'}$ for $d'< g$;
\item $\Theta_{g',k}$ for $g'\leq g-1$ and $k\leq 2g-4$;
\item $\Theta_{g,k}$ for $k\leq 2g-2$.
\end{enumerate}
\end{theo}

{\black Recall that $\Theta_{g,k}$ is the genus $g$ FJRW invariants of insertions $-\frac{2}{5}$ and $\Theta_{g,k}$ may be non-zero only when $k+2-2g\equiv 0 (5)$.
We can see that when $g=2$ only  $\Theta_{2,2}$ is needed,  and when $g=3$ only $\Theta_{3,4}$ is needed.  }

{\black \begin{rema}  As we know,  on using mathematical induction, upon  more numerical datum  the induction   is,  the less effective the computation will be.  We can see from the Theorem that MSP induction for GW invariants is carried out on two numbers, genus and the degree only. Thus this provides a rather effective way to facilitate the induction procedure.

 \end{rema}
}

{\black We can also use the vanishings \eqref{virtual-local} to  get relations among FJRW invariants.}

\begin{theo}[\cite{CLLL2}]\label{nomarking}
Letting $d=(0,d_\infty)$, the vanishings (\ref{virtual-local}) provide relations  
among FJRW invariants $\Theta_{g,k}$.
\end{theo}

{\black These relations are effective in calculating FJRW invariants. For example, for the case of genus $2$, $\{\Theta_{2,k}\}_k$ can be inductively derived from only two unknowns $\Theta_{2,2}$ and $\Theta_{2,7}$.}
\begin{exam}
Computations of $N_{1,1 }$ and $N_{1, 2}$.
\end{exam}
In the Examples \ref{g1d1} and \ref{g1d2},  we listed all the graph types of fixed locus. Using the formulae for $e(N_{F_\Gamma})$ and $[F_\Gamma]^{vir}_{loc}$ in \cite{CLLL2}, we can calculate every term in the summation in (\ref{virtual-local}).

For the genus $1$  degree $1$ case in Example \ref{g1d1}, the contributions from four graph types are { (here $G_i$ is 
the contribution from the graph $\Gamma_i$ in Figure 1)}:  
$$
{ G_1 } = -N_{1,1}, \quad 
{ G_2 } = \frac{9625}{6}, \quad 
{ G_3 } = \frac{-4087}{12}, \quad 
{ G_4 } = -1024.
$$ 
From the equation  (\ref{virtual-local}), the sum of these four numbers should be zero. Thus we obtain $N_{1,1}=\frac{2875}{12}$ which agrees with the known result. 
   
   For the genus $1$ degree $2$ case in Example \ref{g1d2}, the contributions from 15 graph types are { (here
$G_i$ is the contribution from the graph $\Gamma_i$ in Figure 2)}: 
\begin{eqnarray*}
&& { G_1}=N_{1,2},\quad { G_2} =\frac{1106875}{6},\quad { G_3} =\frac{1331125}{12},\quad { G_4}=-\frac{5334375}{2}\\
&& { G_5}=\frac{17206775}{12},\quad { G_6}=355000,\quad { G_7}=-\frac{1018850}{3},\quad { G_8}=\frac{6806875}{4}\\
&& { G_9}=-\frac{12896875}{8} ,\quad { G_{10}}=782000,\quad { G_{11}}=\frac{28966400}{3},\quad { G_{12}}=\frac{4048000}{3}\\
&& { G_{13}}=-\frac{9934400}{3} ,\quad { G_{14}}=-\frac{23116864}{3}, \quad { G_{15}}=12288.\\
 \end{eqnarray*}
 From the equation  (\ref{virtual-local}), the sum of the fifteen numbers above being  zero leads to 
  $N_{1,2}=407125/8$. This is the mathematically verified number by Zinger in \cite{Zi}.
   \qed
\subsection{Speculations}

 Let us look at Theorem \ref{thm-induction} from a different aspect. Inductively we may suppose all GW/FJRW invariants for genus less than $g$ are known. Then for genus $g$, Theorem \ref{thm-induction} reduces the problem of determining  the infinitely  many  GW invariants  $\{N_{g,d}\}_{d=1}^\infty$ to two finite sets of initial datum
 $$ \{N_{g,1},\cdots,N_{g,g-1} \}  \and  \{\Theta_{g,k}\}_{k\leq 2g-2}.$$


We formulate the following speculation:

\medskip

\textit{By suitable choice of positive $d_0$ and $d_\infty$, the relations (\ref{virtual-local}) provide an effective algorithm to determine the first set of initial data $ \{N_{g,1},\cdots,N_{g,g-1} \} $.}

\medskip

 If this is true,  then one is left to determine the second set of initial data $  \{\Theta_{g,k}\}_{k\leq 2g-2}$.  We propose another conjecture about fully determining all FJRW invariants for the quintic,

\begin{conj}\label{conj-2}
The  equations  (\ref{virtual-local})  using $d_0=0$ 
and nonempty $\gamma$'s (i.e. with markings)  give relations  that, together with Theorem \ref{nomarking}, effectively evaluate all $\Theta_{g,k}$.
\end{conj}

 We have verified this conjecture for the case  $\Theta_{2,2}$. Recall that for the case of  genus $2$, this is the only undetermined invariant in Theorem \ref{thm-induction}.

\subsection{Other approaches}
{
The other approach to Witten's proposal is the recent work of Fan, Jarvis and Ruan \cite{FJR3}.
They worked on more general context of gauged linear sigma model, where more general groups $G$ were involved. 
In \cite[Ex.4.2.23]{FJR3}, they took $G=\mathbb C^*\times \mathbb C^*$ and combined the quasi-map technique 
with the P-fields theory to set up the moduli space. 
Incidentally, a closed point of the moduli also consists of a pointed twisted curve $\Sigma\sub\sC$,
two line bundles $\sL$ and $\sN$, and a collection of sections. 
{We point out that despite the similarity, the approach of \cite{FJR3} is different
from ours.} 

The theory in \cite{FJR3} uses  the concept of $\epsilon$-stability, dependent on the real parameter $\epsilon$,
similar to the case of stable quotient \cite{StableQ}. 
The moduli for $\epsilon=0^+$ was constructed in \cite{FJR3}; the case for GW-theory is when $\epsilon=+\infty$,
which is yet to be constructed. 
For  $\epsilon=0^+$ moduli space, the stability (on a point $(\sC,\sL,\sN,\cdots)$) 
requires that $\sL^{-e_1}\otimes \sN^{-e_2}$ is ample on those components of $\sC$ 
for which $\omega^{\log}_\sC$ has degree zero, where $0<{ e_1} <e_2$. 
Coming back to  Example \ref{g1d2}, we see that 
in { $\Gamma_3$} of Figure 2, for the edge $E$ connecting a genus $1$ curve with a genus zero curve,  we have  
$\deg(\omega^{\log}_E)=0$, $\sN|_E\cong \sO$ and $\deg(\sL|_E)=1$. 
Thus $\deg\big((\sL^{-e_1}\otimes \sN^{-e_2})|_E\big)<0$. So curves with the graph type { $\Gamma_3$} in Figure 2 
{ will not} be in  $\epsilon=0^+$ moduli space of \cite[Ex.4.2.23]{FJR3}. 

The $\theta$-parameter in \cite{FJR3} may resemble the $r$-parameter in Witten's vision of CY/LQ
correspondence.} {\black  In out approach, we introduced the new field $\nu=[\nu_1, \nu_2]$ in order to 
``quantize" the Witten's parameter in his phase transition between Calabi-Yau and Landau-Ginzberg 
theories.  We believe that MSP field theory will provide a mathematical theory to realize the vision of Witten.
 We hope that both approaches will be useful for eventual understanding of CY/LG correspondence  in realizing Witten's vision that 
``along a suitable path, there may well be a sharply defined phase transition," 
}
%

Another approach is by Choi and Kiem. In \cite{ChK}, they introduced the moduli spaces of $\epsilon$-stable quasi-maps with P-fields similar to \cite{FJR3} and introduced additional $\delta$-stability to make each wall-crossing more manageable.



\end{document}